\documentclass{article}

\usepackage[latin1]{inputenc}

\usepackage{amsthm,amssymb, amsmath}

\usepackage[all]{xy}

\usepackage{graphicx}

\usepackage[small]{caption}

\usepackage{setspace}

\newcommand{\textaut}{\textsc}

\newcommand{\texttit}{\textit}

\newcommand{\textresto}{\textup}

\newcommand{\matr}[4]{
\left( \begin{array}{cc} #1 & #2 \\ #3 & #4 \end{array} \right)}

\newtheorem{theorem}{Theorem}[section]

\newtheorem{definition}{Definition}[section]

\newtheorem{lemma}[theorem]{Lemma}

\newtheorem{proposition}[theorem]{Proposition}

\newtheorem{corollary}[theorem]{Corollary}

\newtheorem{conjecture}[theorem]{Conjecture}

\newtheorem*{remark}{Remark}

\newcommand{\NN}{\mathbb{N}}

\title{A canonical thickening of $\mathbb{Q}$ and the dynamics of continued fraction transformations}

\author{Carlo Carminati, Giulio Tiozzo}

\begin{document}
\maketitle
% \doublespacing

\begin{abstract}
We construct a countable family of open intervals contained in (0,1]
 whose endpoints are quadratic surds
and such that their union is a full measure set.
We then show that these intervals are precisely 
the monotonicity intervals of the entropy of 
$\alpha$-continued fractions, thus 
proving a conjecture of Nakada and Natsui. 
\end{abstract}

\section{Introduction}

In many areas of mathematics, the space of parameters of a family of mathematical objects 
is itself an object of the same type. A well-known example of this phenomenon in dynamics 
 is the Mandelbrot set, whose local geometry reflects the geometry of the Julia set
of the quadratic polynomial corresponding to a given point.

The goal of this paper is to study
a family of dynamical systems known as 
$\alpha$-continued fraction transformations, showing that the intervals in parameter space 
where a stability condition holds can themselves be  described by means of 
regular continued fraction expansions. 
The family $\{T_\alpha\}_{\alpha \in (0, 1]}$ of $\alpha$-continued fraction transformations 
has been defined in \cite{Nakada81}; the most striking feature of this family may be
that the entropy $h(T_\alpha)$ is not a monotone function of the parameter $\alpha$ (see \cite{LuzziMarmi}), and it is not even 
smooth everywhere.

Rather, Nakada and Natsui (\cite{NakadaNatsui}) showed that the entropy
is locally monotone  on intervals $I$ of parameters which satisfy the following matching condition
\begin{equation} \label{postcrit}
\exists N,M \in \mathbb{N}_+ \ : \quad T^{N}_\alpha(\alpha) = T^{M}_\alpha(\alpha-1) \qquad \forall \alpha \in I
\end{equation}
as well as some other technical conditions. Such intervals will be called \emph{matching intervals}, 
and their union will be referred to as the {\it  matching set}.

In \cite{NakadaNatsui}, Nakada and Natsui exhibited three infinite families of matching
intervals, where the entropy is, respectively, increasing, decreasing, and constant. Moreover, they
conjectured:

\begin{conjecture} \label{conj}
The matching set
has full measure in $(0, 1]$ (hence it is dense).
\end{conjecture}

A numerical study of the conjecture has been carried out in \cite{TiProCarMar}: the goal of this 
paper is to prove the existence of the structures numerically observed there, thus proving conjecture \ref{conj}.

The main tool to analyze the matching set will be regular continued fraction expansions; 
in fact, this matching set can be perfectly described without even mentioning the dynamics of
$\alpha$-transformations. Let us briefly explain why.

It is well known that any rational value $r\in \mathbb{Q}$ can be
expressed as a finite continued fraction expansion of either even or
odd length. This fact, usually perceived as a nuisance, will give us
the chance to perform the following ``natural'' construction:

\begin{enumerate}
 \item For any rational number $r \in \mathbb{Q}\cap (0, 1]$ we consider its two regular continued fraction expansions, namely:
$$r = [0; a_1, \dots, a_n] = [0; a_1, \dots, a_{n} - 1, 1] \qquad a_n \geq 2$$
We will associate to any such $r$ the interval $I_r$ whose endpoints are the quadratic surds 
$$[0; \overline{a_1, \dots, a_n}] \qquad [0; \overline{a_1, \dots, a_{n}-1, 1}]$$
Such an $I_r$ will be called the \emph{quadratic interval} generated by $r$. 

\item
 We will consider the union of all quadratic intervals
$$\mathcal{M} := \bigcup_{r \in \mathbb{Q} \cap (0, 1]} I_r$$
\end{enumerate}

The object of section \ref{quad} will be to understand the structure of the
open dense set  $\mathcal{M}$, which can be summarized in the 

\begin{theorem} \label{first}
The set $\mathcal{M}$ has full Lebesgue measure in $(0,1]$, 
but
its complement has Hausdorff dimension $1$.
\end{theorem}

Although the family of quadratic intervals $\{I_r\}_{r\in \mathbb{Q}}$ will have substantial overlapping, there is a subfamily 
that covers  $\mathcal{M}$ exactly. More precisely,
a quadratic interval $I_r$ will be called \emph{maximal} if it is not properly contained in any other 
quadratic interval. It turns out that
every quadratic interval is contained in some maximal one,
and distinct maximal quadratic intervals do not intersect (lemma \ref{strongmil}): thus
 $\mathcal{M}$ is the
disjoint union of this collection of maximal intervals. 
This suggests that $(0,1]\setminus \mathcal{M}$ should have a
  Cantor-like structure; this is only partially true because
  $(0,1]\setminus \mathcal{M}$ is not perfect.
Indeed, the presence of isolated points is a
    consequence of the {\it period-doubling} phenomenon (see subsection \ref{2x}): if
    $r:=[0;a_1,...a_n]\in \mathbb{Q}$ with $n$ odd and $I_r$ is a
    maximal quadratic interval, then $r':=
    [0;a_1,...a_n,a_1,...a_n]<r$ generates $I_{r'}$ which
    is maximal as well, and the quadratic surd
    $\alpha:=[0;\overline{a_1,...a_n}]$ is a common endpoint, which is obviously
    not contained in any quadratic interval.

In the second part of the paper (section \ref{acf})
we prove that that this set $\mathcal{M}$ is closely connected to the 
matching intervals. More precisely we prove

\begin{theorem} \label{second} Let $a \in \mathbb{Q} \cap (0, 1]$
such that $I_a$ is maximal. Then there exist positive integers $N, M$
such that $$T^{N}_\alpha(\alpha) = T^{M}_\alpha(\alpha-1) \qquad
\forall \alpha \in I_a$$
Moreover, the entropy function $\alpha \mapsto h(T_\alpha)$ is monotone on $I_a$.
\end{theorem}

The proof of the theorem relies on the fact that an \emph{algebraic matching condition} stronger than (1)
holds everywhere on $\mathcal{M}$; by theorem \ref{first}, this condition holds for almost every parameter. 

Moreover, the set defined by the algebraic matching condition contains the matching set defined 
by Nakada and Natsui and the difference between them is countable (see appendix), hence they have the same 
measure and conjecture \ref{conj} follows. 

Our method also gives us an explicit control over the combinatorics of matchings:
given any rational number, we are able to determine which maximal interval it belongs to and the 
\emph{matching exponents} $(N, M)$, hence the local behaviour of entropy (constant, increasing or decreasing).  
Conversely, one can use such knowledge to produce families of matching intervals with prescribed properties.

Finally, section \ref{strtec} contains a few technical tools we use
throughout the paper, including
a criterion to compare purely periodic quadratic surds (String Lemma \ref{STvsTS}) 
and an explicit characterization of either of the finite continued fraction expansions which generate a maximal quadratic interval (lemma \ref{SLSA}).

\vspace{0.3 cm}
It is worth noting that the phenomenon we describe is strongly reminiscent 
of the theory of circle maps (see e.g. \cite{Schuster}, chap. 7.2.): 
in that case, around each rational rotation number,
in the parameter space there is a region ('Arnold tongue') where the dynamics 
is still periodic ('mode-locking'), in such a way that on the critical line
the complement of the union of all Arnold tongues has measure zero (even though its Hausdorff dimension
is strictly smaller than $1$, differently from our case \cite{GS}). 

Recently S. Katok and I. Ugarcovici have studied another family of 
transformations, called $(a,b)$-continued fractions,
 which  seem to share various features with the transformation $T_\alpha$
 (see \cite{KU}): it would be worth investigating more closely the connection between 
these systems in order to see whether the two different approaches can lead to a 
deeper understanding of both.

\section*{Acknowledgements}
We wish to thank the CRM ``E. De Giorgi'' in Pisa, where the bulk of this work was developed, and 
S. Marmi for his constant support. G.T. also wishes to thank C. McMullen for many useful discussions.
This research was partially supported by the M.I.U.R. project ``Dynamical Systems and Applications'' (PRIN 2007B3RBEY)
 
\section{Thickening $\mathbb{Q}$}
 \label{quad}

Let $S = (s_1, \dots, s_n)$ be a finite string of positive integers: we will use the notation 

$$[0; S] := [0; s_1, \dots, s_n] = \frac{1}{s_1 + \frac{1}{ \ddots + \frac{1}{s_n}}}$$
 
Moreover, $\overline{S}$ will be the periodic infinite string $SSS...$ and $[0; \overline{S}]$ the quadratic surd 
with purely periodic continued fraction $[0; \overline{s_1, \dots, s_n}]$. The symbol $|S|$ will 
denote the length of the string $S$. We will denote the denominator of the rational number $r$ as
$\textup{den}(r)$.

\subsection{Pseudocenters}

Let us start out by defining a useful tool in our analysis
of intervals defined by continued fractions.

\begin{lemma} Let $J = ]\alpha, \beta[$, $\alpha, \beta \in \mathbb{R}$, $|\alpha-\beta| < 1$. Then
 there exists a unique rational $p/q \in J$ such that $q = \min \{ q' \geq 1 \ : \ p'/q' \in J \}$.
\end{lemma}

\begin{proof}
Let $d: = \min\{ q\geq 1 \ : \ p/q \in J \} $. If $d = 1$ we are done. Let $d >1$ and assume by contradiction that $\frac{c}{d}$ and $\frac{c+1}{d}$, 
both belong to $J$. Then there exists $k \in \mathbb{Z}$ such that $\frac{k}{d-1} < \frac{c}{d} < \frac{c+1}{d} < \frac{k+1}{d-1}$, hence 
$cd - c- 1 < kd < cd- c$, which is a contradiction since $kd$ is an integer.
\end{proof}

\begin{definition}
The number $\frac{p}{q}$ which satisfies the properties of the previous lemma will be called the \emph{pseudocenter} of $J$.
\end{definition}

\begin{lemma} \label{pseudoc}
Let $\alpha, \beta \in ]0, 1[$ be two irrational numbers with c.f. expansions $\beta = [0; S, b_0, b_1, b_2, \dots]$
and $\alpha = [0; S, a_0, a_1, a_2, \dots]$, where $S$ stands for a finite string of positive integers. Assume $b_0 > a_0$. Then the pseudocenter of the interval $J$ with endpoints $\alpha$ and $\beta$ is
$$r = [0; S, a_0 + 1] ( = [0; S, a_0, 1])$$ 
\end{lemma}

\begin{proof}
Suppose there exists $s \in \mathbb{Q} \cap J$ with den($s$) $<$ den ($r$). Since $s \in J$, then $s = [0 ; S, s_0, s_1, \dots, s_k]$
 with $a_0 \leq s_0 \leq b_0$ and $k \geq 0$. The choice $s_0 \geq a_0 + 1$ gives rise to den($s$) $\geq$ den($r$), so $s_0 = a_0$. On the other hand, 
$[0; S, a_0]$ does not belong to the interval, so $k \geq 1$ and $s_1 \geq 1$, still implying den($s$) $\geq$ den($r$).
\end{proof}

\subsection{Quadratic intervals}

\begin{definition}
Let $0< a < 1$ be a rational number with c.f. expansion
$$a = [0; a_1, \dots, a_N] = [0; a_1, \dots, a_N-1, 1], \ a_N \geq 2$$
We define the \emph{quadratic interval} $I_a$ associated to $a$ to be the open interval with endpoints 
\begin{equation}
 \label{endpoints}
[0; \overline{a_1, \dots, a_{N-1}, a_N}] \quad \textup{and} \quad [0; \overline{a_1, \dots, a_{N-1}, a_N-1, 1}]
\end{equation}
Moreover, we define 
$I_1 := (\frac{\sqrt{5}-1}{2}, 1]$ (recall that $\frac{\sqrt{5}-1}{2}= [0; \overline{1}]$).
\end{definition}

Note that the ordering of the endpoints in \eqref{endpoints} depends on the parity of $N$: 
given $a \in \mathbb{Q}$, we will denote by $A^+$ and $A^-$ the two strings of positive integers which represent 
$a$ as a continued fraction, with the convention that $A^+$ is the string of \emph{even} length and $A^-$ the string of \emph{odd} length, 
so that 
$$I_a = ([0; \overline{A^-}], [0; \overline{A^+}]),  \qquad a = [0; A^+] = [0; A^-] $$

{\bf Example} 

If $a = \frac{1}{3} = [0; 3] = [0; 2, 1]$, $[0; \overline{A^+}] = [0; \overline{2, 1}]$, $[0; \overline{A^-}] = [0; \overline{3}]$, $I_a = (\frac{\sqrt{13}-3}{2}, \frac{\sqrt{3}-1}{2})$

\noindent Note that $a$ is the pseudocenter of $I_a$, hence $I_a = I_{a'} \Leftrightarrow a = a'$.  

\begin{lemma} \label{lemmacontained}
\begin{enumerate}
\item If $\xi \in \overline{I}_a$, then $a$ is a convergent to $\xi$. 
\item If $I_a \cap I_b \neq \emptyset$, then either $a$ is a convergent
%\footnote{We will say that a rational $r$ is a convergent to another 
%rational $q$ if it is a convergent to at least one 
%of the two possible c.f. expansions of $q$.} 
to $b$ or $b$ is a convergent to $a$. 
\item If $I_a \subsetneq I_b$ then $b$ is convergent to $a$, hence $den(a) < den(b)$. 
\end{enumerate}
\end{lemma}

\begin{proof}
1. Since $\xi \in I_a$, either $\xi = [0; a_1, \dots, a_N, \dots]$ or $\xi = [0; a_1, \dots, a_N-1, \dots]$. In the first case the claim holds; 
in the second case one has to notice that neither $[0; a_1, \dots, a_N-1]$ nor all elements of the form $[0; a_1, \dots, a_N-1, k, \dots]$
 with $k \geq 2$ belong to $I_a$, so $k = 1$ and $a$ is a convergent of $\xi$.

2. Fix $\xi \in I_a \cap I_b$. By the previous point, both $a$ and $b$ are convergents of $\xi$, hence the rational with 
the shortest expansion is a convergent of the other.

3. From 1. since $a \in I_a \subseteq I_b$.
\end{proof}

\begin{definition}
A quadratic interval $I_a$ is \emph{maximal} if it is not properly contained in any $I_b$ with $b \in \mathbb{Q}\cap ]0, 1]$. 
\end{definition}

The interest in maximal quadratic intervals lies in the 

\begin{proposition} \label{unique}
Every quadratic interval $I_a$ is contained in a unique maximal quadratic interval.
\end{proposition}

A good way to visualize the family of quadratic intervals is to plot, for any rational $a$, 
the geodesic $\gamma_a$ on the hyperbolic upper half plane with the same endpoints as $I_a$, as in the following picture:
one can see the maximal intervals corresponding to the ``heighest'' geodesics, in such a way that 
every $\gamma_a$ has some maximal geodesic (possibly itself) above it and no two maximal $\gamma_a$ intersect.

\begin{center}
\includegraphics[scale = 0.5]{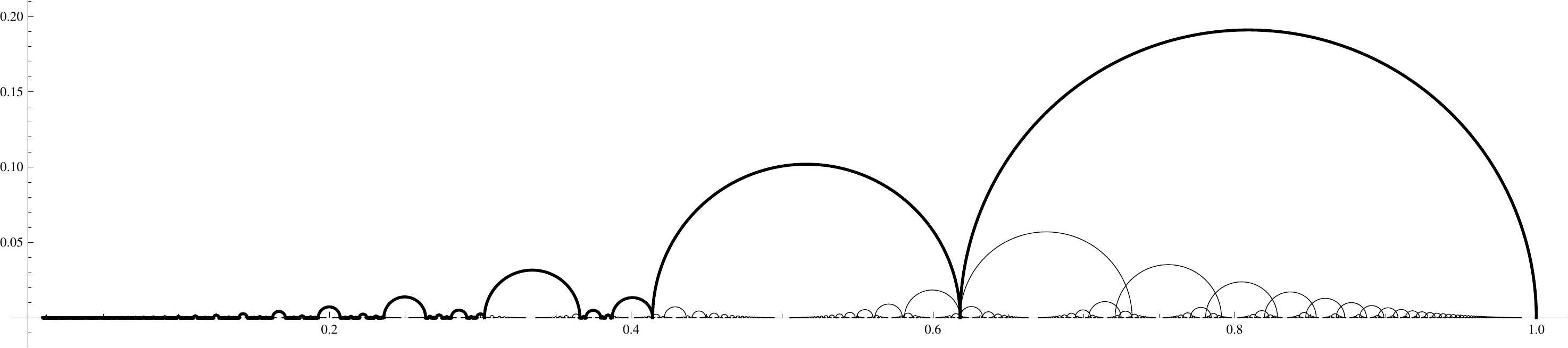}
\end{center}

The proof of proposition \ref{unique} will be given in two lemmas:
\begin{lemma}
Every quadratic interval $I_a$ is contained in some maximal quadratic interval.
\end{lemma}

\begin{proof}
If $I_a$ were not contained in any maximal interval, then there would exist an infinite chain $I_{a} \subsetneq I_{a_1} \subsetneq I_{a_2} \subsetneq \dots$
of proper inclusions, hence by the lemma every $a_i$ is a convergent of $a$, but rational numbers can only have a finite number of convergents.
\end{proof}

\begin{lemma}\label{strongmil}
If $I_a$ is maximal then for all $a'\in \mathbb{Q}\cap (0,1)$ 
$$I_a \cap I_{a'} \neq \emptyset \ \Rightarrow I_{a'}\subset I_a ,$$
and equality holds iff $a=a'$. In particular, distinct maximal intervals do not intersect.
\end{lemma}

\begin{proof}

We need the following lemma, which we will prove in section \ref{strtec}:
\begin{lemma} \label{maxintlemma}
If $I_a \cap I_b \neq \emptyset$, $I_a \setminus I_b \neq \emptyset$ and $I_b \setminus I_a \neq \emptyset$, then 
either $I_a$ or $I_b$ is not maximal. 
\end{lemma}
Let now $I_{a_0}$ be the maximal interval which contains $I_{a'}$. Since $I_a \cap I_{a_0} \neq \emptyset$, 
by lemma \ref{maxintlemma} either $I_a \subseteq I_{a_0}$ or $I_{a_0} \subseteq I_a$, hence by maximality 
$I_a = I_{a_0}$ and $I_{a'} \subseteq I_a$. Since $a$ is the pseudocenter of $I_a$, $I_a = I_{a'} \Rightarrow a = a'$.
\end{proof}
\subsection{Hausdorff dimension}

In this section we prove theorem \ref{first}, which states that the
exceptional set $\mathcal{E}:=]0,1]\setminus \mathcal{M}$ has zero
    Lebesgue measure but Hausdorff dimension equal to 1. The key tool
    of the proof is the following lemma, which establishes a
    connection between $\mathcal{E}$ and numbers of bounded type.
\noindent 
\begin{lemma} \label{infbound}
\begin{enumerate}
\item[(i)] Let $\xi \in \mathcal{E}=]0,1]\setminus \mathcal{M}$. Then $\xi$  is irrational and  $\xi = [0; a_1, \dots, a_n, \dots]$ with $a_j\leq a_1$ for all $j\in \mathbb{N}_+$
 \item[(ii)] Let $\xi = [0; a_1, \dots, a_n, \dots]$ be an irrational number such that $a_k \leq a_1 -1 $ for all $k \geq 2$.
Then $\xi$ does not belong to any $I_a$ for any $a \in \mathbb{Q}\cap(0, 1]$.
\end{enumerate}
\end{lemma}

\begin{proof}
Since $\xi \notin \mathcal{M}$ then $\xi \notin \mathbb{Q}$. If $\xi$ has the infinite c.f. expansion $\xi = [0; a_1, \dots, a_n, \dots]$ with $a_k>a_1$ for some $k \in \mathbb{N}_+$ then $x$ lies between $r:=  [0; a_1, \dots, a_{k-1}]$
 and $\alpha:=[0; \overline{a_1, \dots, a_{k-1}}]$; therefore $x\in I_r\subset \mathcal{M}$. 
Let $a = [0; A^+] = [0; A^-]$, so that $I_a = ([0; \overline{A^-}], [0; \overline{A^+}])$.
If $\xi \in I_a$, by lemma \ref{lemmacontained} $a$ is a convergent of $\xi$, so either
 $$\xi = [0; A^+, \dots] \quad \textup{ or}\quad \xi = [0; A^-, \dots]$$
In the first case $\xi = [0; A^+, s, \dots]$ with $s < a_1$, so 
$\xi > [0; \overline{A^+}] = [0; A^+, a_1, \dots]$; in the second one, 
$\xi = [0; A^-, s, \dots]$ with $s < a_1$ and therefore
$\xi < [0; \overline{A^-}] = [0; A^-, a_1, \dots]$. 
\end{proof}

\begin{proof} 
(of Theorem \ref{first})
Lemma \ref{infbound} implies that $\mathcal{E}$ is contained in the set of numbers of bounded type, hence it has Lebesgue measure zero.
  
On the other hand, let $N \geq 1$, and define  
$$C_{N} := \{x = [0; a_1, \dots] \ \mid \ a_k \leq N \ \forall k \geq 1\}$$
$$E_N := \left[\frac{1}{N+1}, \frac{1}{N}\right) \cap \mathcal {E}$$
By lemma \ref{lemmanbounded} and lemma \ref{infbound}
$E_N \subseteq C_N$ 
and by lemma \ref{infbound}, for $N \geq 2$,
$E_N \supseteq \phi(C_{N-1})$
where $\phi(x) := x \mapsto \frac{1}{N+x}$. 
Since $\phi$ is a bi-Lipschitz map, it preserves Hausdorff dimension, so 
$$\textup{dim}_H {C_{N-1}} = \textup{dim}_H {\phi(C_{N-1})} \leq \textup{dim}_H E_N \leq \textup{dim}_H C_N$$
Since it is well-known (\cite{Jarnik}) that $\sup_{N \rightarrow \infty} \textup{dim}_H \ C_N = 1$ and $\mathcal{E}=\cup_N E_N$, 
the claim follows.
\end{proof}

\begin{remark}
A similar way of stating the same result would be to say that for every $\frac{p}{q} \in \mathbb{Q}\cap \left(\frac{1}{N+1}, \frac{1}{N}\right)$ 
$$B\left( \frac{p}{q}, \frac{1}{(N+2)q^2} \right) \subseteq I_{p/q} \subseteq B\left( \frac{p}{q}, \frac{1}{(N-1)q^2} \right)$$
This means that in any fixed subinterval $(\frac{1}{N+1}, \frac{1}{N})$ the size of the geodesic over $I_{p/q}$
is comparable to the diameter of the horocycles
$\partial B(\frac{p}{q} + \frac{\imath}{Nq^2}, \frac{1}{Nq^2})$ (which, for any fixed $N$, all lie in the same $SL_2(\mathbb{Z})$-orbit).
The picture shows this comparison for $N = 10$.
\end{remark}

\begin{center}
\includegraphics[scale=0.25]{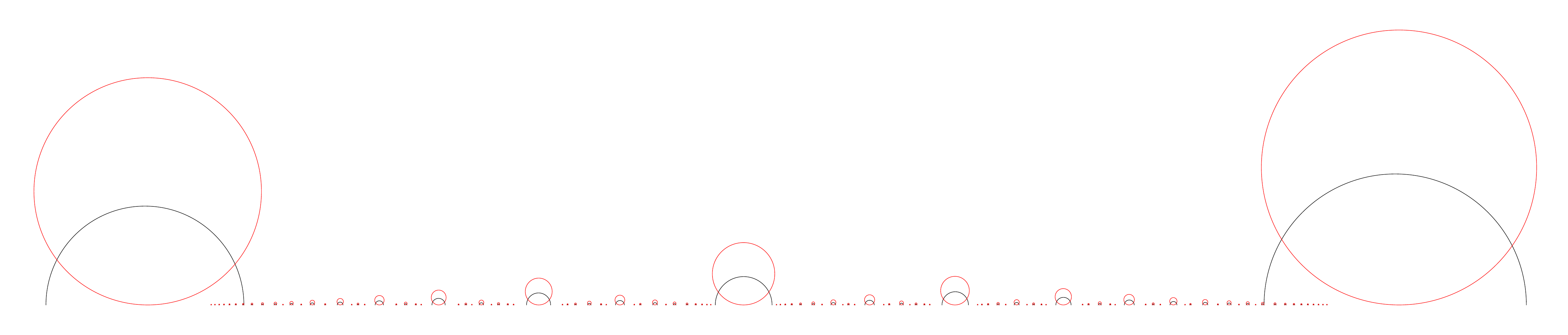}
\end{center}

\subsection{The bisection algorithm}

We will now describe an algorithmic way to produce all maximal intervals, as announced in \cite{TiProCarMar}, sect. 4.1. 
This will also provide an alternative proof of the fact the $\mathcal{M}$ has full measure.

Let $\mathcal{F}$ be a family of disjoint open
intervals which \emph{accumulate only at }$0$, i.e. such that
for every $\epsilon > 0$ the set $\{J \in \mathcal{F} : J\cap[\epsilon, 1] \neq \emptyset \}$ is finite, 
and denote $F = \cup_{J \in \mathcal{F}} J$. 
The complement $]0, 1] \setminus F$ will then be a countable union of closed disjoint intervals $C_j$, 
which we refer to as \emph{gaps}. Note that some $C_j$ may well be a single point.
To any gap which is not a single point we can associate its pseudocenter $c\in \mathbb{Q}$ as defined in the previous sections, 
and moreover consider the interval $I_c$ associated to this rational value. The following proposition applies.

\begin{proposition}
Let $I_a$ and $I_b$ be two maximal intervals such that the gap between them is not a single point, and let $c$ be the 
pseudocenter of the gap. Then $I_c$ is a maximal interval and it is disjoint from both $I_a$ and $I_b$.
\end{proposition}

\begin{proof}
Pick $I_{c_0}$ maximal such that $I_c \subseteq I_{c_0}$, so by lemma  \ref{lemmacontained} 
  $den(c_0) \leq den(c)$. On the other hand, since maximal intervals do not intersect, 
then $I_{c_0}$ is contained in the gap and since $c$ is pseudocenter, then $den(c) \leq den(c_0)$ and equality
 holds only if $c = c_0$. 
\end{proof}

The proposition implies that if we add to the family of maximal intervals $\mathcal{F}$ all intervals which arise as gaps between adjacent 
intervals then we will get another family of maximal (hence disjoint) intervals, and we can iterate the procedure.

For instance, let us start with the collection $\mathcal{F}_1 := \{ I_{1/n}, n\geq 1\}$. All these intervals are maximal, 
since the continued fraction of their pseudocenters has only one digit (apply lemma \ref{lemmacontained}). 

Let us construct the families of intervals $\mathcal{F}_n$ recursively as follows:
$$\mathcal{G}_n := \{C \textup{ connected component of }]0,1]\setminus F_n\}$$

$$\mathcal{F}_{n+1} := \mathcal{F}_n \cup \{ I_r \ : \  r \textup{ pseudocenter of }C, C\in \mathcal{G}_n, C \textup{ not a single point }  \}$$
(where $F_n$ denotes the union of all intervals belonging to $\mathcal{F}_n$).

It is thus clear that the union $\mathcal{F}_\infty := \bigcup \mathcal{F}_n$ will be a countable family of maximal intervals. 
The union of all elements of $\mathcal{F}_\infty$ will be denoted by $F_\infty$; its complement (the set of numbers which do not 
belong to any of the intervals produced by the algorithm) has the following property:

\begin{lemma} \label{lemmanbounded}
$]0,1[\setminus F_\infty$ consists of irrational numbers of bounded type; more precisely,
the elements of $(\frac{1}{n+1}, \frac{1}{n}]\setminus F_\infty$ have partial quotients bounded by $n$.
\end{lemma}

\begin{proof}
Let $\gamma=[0;c_1, c_2,..., c_n, ...]\notin F_\infty$; we claim that $c_k\leq c_1$ for
 all $k \in \mathbb{N}$. Since $\gamma \notin F_\infty$, $\forall n \geq 1$ we can choose $J_n \in \mathcal{G}_n$
such that $\gamma \in J_n$. Clearly, $J_{n+1} \subseteq J_n$.
Furthermore, $\gamma$ cannot be contained in either $I_{\frac{1}{c_1}}$ nor $I_{\frac{1}{c_1+1}}$, so all $J_n$ are produced 
by successive bisection of the gap $([0; \overline{c_1, 1}], [0; \overline{c_1}])$, hence by lemma \ref{pseudoc}
for every $n$, the endpoints of $J_n$ are quadratic surds with c.f. expansion bounded by $c_1$.
It may happen that there
 exists $n_0$ such that $J_n=\{\gamma\} \ \forall n\geq n_0$, so $\gamma$ is an endpoint of 
$J_{n_0}$, hence it is irrational and $c_1$-bounded. Otherwise,
let $p_n/q_n$ be the pseudocenter of $J_n$; by uniqueness of the pseudocenter, 
$\textup{diam }J_n \leq 2/q_n$, and $q_{n+1} > q_n$ since $J_{n+1} \subseteq J_n$.
This implies $\gamma$ cannot be rational, since the minimum denominator of a rational sitting
 in $J_n $ is $q_n \to +\infty$. Moreover, $\textup{diam }J_n \rightarrow 0$, 
so $\gamma$ is limit point of endpoints of the $J_n$, which are $c_1$-bounded, hence 
$\gamma$ is also $c_1$-bounded.
\end{proof}

\begin{proposition}\label{maximalbybisection}
The family $\mathcal{F}_\infty$ is precisely the family of all maximal intervals; hence $F_\infty = \mathcal{M}$.
\end{proposition}

\begin{proof}
If $I_c$ a maximal interval does not belong to $\mathcal{F}_\infty$, then its pseudocenter belongs to the 
complement of $F_\infty$, but the previous lemma asserts that this set does not contain any rational. 
\end{proof}

Note that proposition \ref{maximalbybisection} and lemma \ref{lemmanbounded} provide another way of seeing that 
the complement of $\mathcal{M}$ consists of numbers of bounded type, hence it has full measure.

\subsection{Maximal intervals and strings} \label{coding}

In order to get a finer control on the maximality properties of quadratic intervals, we introduce 
a systematic description of the continued fraction expansions in terms of strings and develop a few 
tools in order to characterize the expansions of those rational numbers which give rise to maximal intervals.

Let us start out with some notation. If $S = (s_1, \dots, s_n)$ is a finite string of positive integers 
and $x$ a real number, we will denote

$$[0; S] := \frac{1}{s_1 + \frac{1}{ \ddots + \frac{1}{s_n}}} \qquad 
[0; S + x] := \frac{1}{s_1 + \frac{1}{ \ddots + \frac{1}{s_n+x}}}$$
 
We will also introduce a total ordering on the space of finite strings of given length:
given two distinct finite strings $S$ and $T$ of equal length, let $l := \min \{i : S_i \neq T_i \}$.We will set
$$S < T := \left\{\begin{array}{l} S_l < T_l \textup{ if }l \equiv 0 \mod 2 \\
                     S_l > T_l \textup{ if }l \equiv 1 \mod 2 
                    \end{array} \right.$$
The exact same definition also gives a total ordering on the space of infinite strings. Note that if 
$S$ and $T$ have equal length $L \in \mathbb{N} \cup \{\infty\}$, 
$$S < T \Leftrightarrow [0; S] < [0; T]$$
i.e. this ordering can be obtained by pulling back the order structure on $\mathbb{R}$, via
 identification of a string with the value of the corresponding c.f.

The following lemma is the essential tool used to compare two purely periodic infinite strings:

\begin{lemma}\label{STvsTS}
Let $S, \ T$ be two nonempty, finite strings. Then the pair of infinite strings
$\overline{S}$, $\overline{ T}$ is ordered in the same way as the pair $ST$, $TS$; namely 
 $$ ST \gtreqqless TS \iff \overline{S} \gtreqqless \overline{ T}.$$
\end{lemma}

Finally, we can give an explicit characterization of the c.f. expansion of those rationals 
which are pseudocenters of maximal intervals:

\begin{proposition} \label{SLSA}
Let $a = [0; A] \in \mathbb{Q} \cap (0, 1]$. The following are equivalent:
\begin{itemize}
 \item[(i)] $I_a$ is maximal.
\item[(ii)] If $A = ST$ with $S$, $T$ finite nonempty strings, then either
$ST < TS$ or $ST = TS$ with $T = S$, $|S|$ odd
\end{itemize}
  
Moreover, if $[0; ST]$ is maximal, then $[0; T] > [0; ST]$
\end{proposition}

For the sake of readability, we postpone the proofs of these results to section \ref{strtec}.

\section{Application to $\alpha$-continued fractions} \label{acf}

After having investigated the properties of the maximal set itself, this 
section will be devoted to studying its relation with the parameter space of 
$\alpha$-continued fractions.

\subsection{Matching intervals} \label{matchint}

Let $\alpha \in (0, 1]$. Recall that the $\alpha$-continued fraction expansion is given by the map
$T_\alpha : [\alpha-1, \alpha] \rightarrow [\alpha-1, \alpha]$ defined by $T_\alpha(0)=0$ and
$$T_{\alpha}(x) = \frac{\epsilon_{\alpha}(x)}{x} - c_{\alpha}(x) \qquad \textup{for }x \neq 0$$ 
with 
$$\epsilon_{\alpha}(x) := \textup{Sign}(x) \qquad c_{\alpha}(x) := \left\lfloor \frac{1}{|x|} + 1 - \alpha \right\rfloor$$
Moreover, one can represent the encoding with the matrices in $GL(2, \mathbb{Z})$

$$M_{\alpha, x, n} = \matr{0}{\epsilon_{\alpha}(x)}{1}{c_{\alpha}(x)} \dots \matr{0}{\epsilon_{\alpha}(T_\alpha^{n-1}(x))}{1}{c_{\alpha}(T_\alpha^{n-1}(x))} = \matr{p_{n-1, \alpha}(x)}{p_{n,\alpha}(x)}{q_{n-1, \alpha}(x)}{q_{n,\alpha}(x)}$$
so that
\begin{equation} \label{invert}
 x = \frac{p_{n-1, \alpha}(x) x_n + p_{n, \alpha}(x)}{q_{n-1, \alpha}(x) x_n + q_{n,\alpha}(x)} \quad \textup{with }x_n = T^n_{\alpha}(x).
\end{equation}

We will be interested in the metric entropy $h(T_\alpha)$ of these transformations as a function of $\alpha$;
in \cite{NakadaNatsui}, a series of \emph{matching conditions} were introduced in order to 
define intervals in the parameter space where the entropy function $\alpha \mapsto h(T_\alpha)$ is monotone. 
In the same spirit, we will define

\begin{definition}
The value $\alpha \in ]0, 1]$ is said to satisfy an \emph{algebraic matching condition} of order $(N, M)$
when the following matrix identity holds:
\begin{equation*}
M_{\alpha, \alpha, N} = \matr{1}{1}{0}{1} M_{\alpha, \alpha-1, M} \matr{-1}{0}{1}{1} \qquad (N,M)_{\rm alg} 
\end{equation*} 
\end{definition}

\noindent We will be interested in the set
$$\mathcal{M}_{\rm alg} = \{ \alpha \in ]0, 1] \textup{ s.t. } \exists N, M
    \in \mathbb{N} \ : \alpha \textup{ satisfies } (N,M)_{\rm alg} \}$$

To get some intuition of what this condition means from a dynamic point of view, 
one should note that $(N,M)_{\rm alg}$ implies
$$T_\alpha^{N+1}(\alpha) = T_\alpha^{M+1}(\alpha-1)$$
The formal proof of this result is given in the appendix, together with a thorough discussion
of the relationship between our algebraic matching condition and the conditions originally considered by Nakada 
and Natsui. 

The main result will be:

\begin{theorem} \label{main}
Let $a \in \mathbb{Q} \cap (0, 1]$ such that $I_a$ is maximal, and let $a = [0; a_1, \dots, a_n]$, $n$ even.
If we define 
$$N := \sum_{j \ \textup{even}}a_j  \qquad M := \sum_{j \ \textup{odd}} a_j$$
then for every $x \in I_a$, the matching condition $(N,M)_{\rm alg}$ holds.
\end{theorem}

\begin{corollary}
$\mathcal{M}_{\rm alg}$ has full Lebesgue measure in $(0, 1]$. 
\end{corollary}

\begin{proof}
By theorem \ref{main}, $\mathcal{M}_{\rm alg}$ contains $\mathcal{M}$, which has full measure by theorem \ref{first}.
\end{proof}

Since it can be proved (see appendix) that the difference between $\mathcal{M}_{\rm alg}$ and the matching set defined by 
Nakada and Natsui is countable, this also establishes conjecture \ref{conj}.

\subsection{Anatomy of maximal orbits}

The first step in the proof of theorem \ref{main} will be to describe explicitly the first few steps 
of the orbit of any point inside a maximal interval $I_a$: we will start by establishing the 

\begin{lemma} \label{anatomy}
Let $a \in \mathbb{Q} \cap (0, 1]$ be  the pseudocenter of a maximal $I_a = (\alpha^-, \alpha^+)$. 
\begin{enumerate}
 \item 
Let $a \leq x < \alpha^+$, so that we can write $x = [0; a_1, \dots, a_n + y]$ with $0 \leq y < \alpha^+$, 
$a = [0; a_1, \dots, a_n]$ with $n \equiv 0 \mod 2$. Then 
$$[-1; b, a_{k+1}, \dots, a_n + y] > \alpha^+ -1 \quad \forall 1 \leq b \leq a_k, \ 1 < k \leq n$$
 \item
Let $\alpha^- < x \leq a$, so that $x = [0; a_1, \dots, a_n + y]$ with $0 \leq y < \alpha^-$, 
$a = [0; a_1, \dots, a_n]$ with $n \equiv 1 \mod 2$ (note this is the representation of $a$ in c.f. \emph{other} than 
the one given in the previous point). Then 
  
$$[-1; b, a_k, \dots, a_n + y] > \alpha^+ - 1 \quad \forall 1 \leq b \leq a_k, \ 1 < k \leq n$$
\end{enumerate}
\end{lemma}

\begin{proof}
1. Let $S := (a_1, \dots, a_{k-1})$, $T := (a_k, \dots, a_n)$ and $c := [0; T]$.

By lemma \ref{SLSA} and \ref{STvsTS},
$$TS \geq ST \Rightarrow \overline{TS} \geq \overline{ST} \Rightarrow [0; \overline{TS}] \geq [0; \overline{ST}]$$ 
Moreover,
$$ TS \geq ST \Rightarrow TST \geq STT \Rightarrow \overline{T} \geq \overline{ST} $$
Now, $I_c \cap I_a = \emptyset$ since $I_a$ is maximal and the denominator of $c$ 
is smaller than the denominator of $a$, hence $[0; T] > [0; \overline{ST}]$.
Since $b \leq a_k$ and $0 \leq y < \alpha^+$, for $k$ even we have 
$$[-1; b, a_{k+1}, \dots, a_n + y ] \geq [-1; T, y] > [-1; T, \alpha^+] = [-1; \overline{TS}] \geq [-1; \overline{ST}] = \alpha^+ -1$$
and for $k$ odd,
$$[-1; b, a_{k+1}, \dots, a_n + y ] \geq [-1; T, y] \geq [-1; T] > [-1; \overline{ST}] = \alpha^+ -1$$

2. Let $S := (a_1, \dots, a_{k-1})$, $T := (a_k, \dots, a_n)$ and $c := [0; T]$.
If $k$ is odd
$$TS \geq ST \Rightarrow TTS \leq TST \Rightarrow \overline{T} \leq \overline{TS}$$
Moreover, $\overline{T} \geq \overline{ST}$ as in the previous point, and since $I_a \cap I_c = \emptyset$, 
then $[0; \overline{T}] \geq \alpha^+$, so $[0; \overline{TS}] \geq [0; \overline{T}] \geq \alpha^+$; hence,
$$[-1; b, a_{k+1}, \dots, a_n + y] \geq [-1; T, y] > [-1; T, \alpha^-] = [-1; \overline{TS}] \geq \alpha^+ -1.$$
For $k$ even, by the last point of proposition \ref{SLSA}, $[0; T] > [0; ST]$, and since $I_a \cap I_c = \emptyset$, 
$[0; T] > \alpha^+$; thus,
$$ [-1; b, a_{k+1}, \dots, a_n + y] \geq [-1; T, y] \geq [-1; T] >  \alpha^+ -1$$

\end{proof}

An immediate corollary is the explicit description of the orbit 
of the pseudocenter which explains an empirical 
rule given in \cite{TiProCarMar}.

\begin{corollary}\label{orbitpseudo}
Let $a:=[0;a_1,a_2,...a_n]$, $(n\geq 1)$ and let $I_a$ be maximal; then the orbits of $a$ and $a-1$ are as follows:
$$
\begin{array}{lcl}
a=[0;a_1,a_2,...a_n] & \ & a-1=[-1;a_1,a_2,...a_n]\\
T_a(a)=[-1;a_2,...a_n] &\ & T_a(a-1)=[-1;a_1-1,a_2,...a_n]\\
... &\ & ...\\
T_a^{a_2}(a)=[-1;1,a_3,...a_n] &\ & T_a^{a_1-1}(a-1)=[-1;1,a_2,...a_n]\\
T_a^{a_2+1}(a)=[-1;a_4,...a_n] & \ & T_a^{a_1}(a-1)=[-1;a_3,...a_n]\\
... &\ & ...\\
T_a^{N}(a)=0 &\ & T_a^{M}(a-1)=0\\
\end{array}
$$
where (see also \cite{TiProCarMar}, pg. 23)
$$
\begin{array}{lll}
\displaystyle N=\sum_{j \mbox{ even}}a_j, &\displaystyle  M=\sum_{j \mbox{ odd}}a_j, & \mbox{ if n is even}\\

\displaystyle N=1+\sum_{j \mbox{ even}}a_j, &\displaystyle  M=-1+\sum_{j \mbox{ odd}}a_j, & \mbox{ if n is odd}\\[16pt]

\end{array}$$
\end{corollary}

We will now prove that an algebraic matching condition holds for any pseudocenter of a maximal interval.

\begin{proposition} \label{matchrat}
Let $a \in \mathbb{Q} \cap (0, 1]$ so that $I_a$ is maximal, 
and let $N$ and $M$ be given by the previous corollary. Then $a$ satisfies the matching condition $(N,M)_{\rm alg}$. 
\end{proposition}

\begin{proof}
We will make use of the following lemma:

\begin{lemma} \label{incrden}
For $\alpha < \frac{\sqrt{5}-1}{2}$, one has $q_{n+1, \alpha}(x) > q_{n, \alpha}(x) \geq 1$ 
for every $n \geq 0$ and every $x \in [\alpha-1, \alpha]$.
\end{lemma}

\begin{proof}
By definition, $q_{0, \alpha}(x) = 1$ and $q_{1, \alpha}(x) = c_{1, \alpha}(x) \geq 2 $ (the latter only for $\alpha < \frac{\sqrt{5}-1}{2}$). By induction, 
$q_{n+1, \alpha}(x) = c_{n+1, \alpha}(x) q_{n,\alpha}(x) + \epsilon_{n+1, \alpha}(x) q_{n-1, \alpha}(x) \geq 2 q_{n, \alpha}(x) - q_{n-1, \alpha}(x) > q_{n, \alpha}(x)$.
\end{proof}

Since it is easy to see that all values of $\alpha > \frac{\sqrt{5}-1}{2}$ satisfy a matching condition of order 
$(1, 2)$, we can restrict our attention to the case in which we can apply lemma \ref{incrden}. 
We will denote $p_{k} := p_{k, \alpha}(\alpha)$ and $p'_{k} := p_{k, \alpha}(\alpha-1)$. 
Let $(N, M)$ be given by corollary \ref{orbitpseudo}, such that 
$$T_{a}^N(a) = 0 \ \textup{ and } \ T_{a}^M(a-1) = 0$$ 
By equation  \eqref{invert}, 
$$a = p_{N}/q_{N} \qquad a-1 = p'_{M}/q'_{M}$$
and since $gcd(p_{N}, q_{N}) = gcd(p'_{M}, q'_{M}) = 1$ (because $\det M_{a, x, k} = \pm 1$),
\begin{equation} \label{qid1} 
q_{N} = q'_{M} \qquad p_{N} = p'_{M} + q'_{M}
\end{equation}

Now, corollary \ref{orbitpseudo} implies $\epsilon_{a}(T_a^i(a)) = \epsilon_{a}(T_a^{j}(a-1)) = -1$ for $1 \leq i \leq N-1$, $1 \leq j \leq M-1$, hence 
$$\det M_{a, a, N} = -1 \qquad \det M_{a,a-1, M}= 1$$
by writing out the two determinants and summing up 
$$p_{N-1}q_N-p_Nq_{N-1} + p'_{M-1}q'_M - p'_M q'_{M-1} = 0$$
and by using \eqref{qid1}
$$q'_M(p_{N-1} + p'_{M-1} - q_{N-1}) = p'_M(q'_{M-1} + q_{N-1})$$
Now, $q'_M$ and $p'_M$ are coprime, hence $q'_M | (q'_{M-1} + q_{N-1})$, 
and by lemma \ref{incrden}, $0 < q'_{M-1} + q_{N-1} < 2 q'_M$, so 
$$q'_M = q'_{M-1} + q'_{N-1} \qquad p'_M = p_{N-1} + p'_{M-1} - q_{N-1}$$ 
which yields precisely the algebraic matching condition
$$M_{a,a, N} =  \matr{1}{1}{0}{1} M_{a, a-1, M} \matr{-1}{0}{1}{1}$$
\end{proof}

The final step will be to prove that all points in $I_a$ have the same convergents as the pseudocenter.

\begin{lemma} \label{sameorbit}
Let $I_a$ be maximal, and $x \in I_a$, $N$, $M$ as in corollary \ref{orbitpseudo}. Then
$$\begin{array}{ll} M_{x, x, k} = M_{a, a, k} & \forall 1 \leq k \leq N \\
   M_{x, x-1, h} = M_{a, a-1, h} & \forall 1 \leq h \leq M
\end{array}$$
\end{lemma}

\begin{proof}
If $x \geq a$, we can write $x = [0; A + y]$ with $|A| \equiv 0 \mod 2$ and $0 \leq y < \alpha^+$; 
from corollary \ref{orbitpseudo} 
$$
\begin{array}{lcl}
x=[0;a_1,a_2,...a_n + y] & \ & x-1=[-1;a_1,a_2,...a_n + y]\\
M_{a, a, 1}^{-1}(x)=[-1;a_2,...a_n + y] &\ & M_{a, a-1, 1}^{-1}(x-1)=[-1;a_1-1,a_2,...a_n + y]\\
... &\ & ...\\
M_{a, a, a_2}^{-1}(x)=[-1;1,a_3,...a_n+ y] &\ & M_{a, a-1, a_1-1}^{-1}(x-1)=[-1;1,a_2,...a_n+y]\\
M_{a, a, a_2+1}^{-1}(x)=[-1;a_4,...a_n+ y] & \ & M_{a, a-1, a_1}^{-1}(x-1)=[-1;a_3,...a_n+y]\\
... &\ & ...\\
M_{a, a, N}^{-1}(x)= [-1; 1+ y] &\ & M_{a, a-1, M}^{-1}(x-1)= y\\
\end{array}
$$
and again from the lemma, 
$$M_{a, a, k}^{-1}(x) \in (\alpha^+-1, 0) \subseteq (x-1, 0) \qquad 1 \leq k \leq N$$ 
$$M_{a, a-1, h}^{-1}(x) \in (\alpha^+-1, 0) \subseteq (x-1, 0) \qquad 1 \leq h \leq M-1$$
hence
$$\begin{array}{ll} M_{x, x, k} = M_{a, a, k} & 1 \leq k \leq N \\
   M_{x, x-1, h} = M_{a, a-1, h} & 1 \leq h \leq M-1.
  \end{array}$$
To prove the claim we are left with considering 
$$ M_{a, a-1, M}^{-1}(x-1) = y$$
Since
$$0 < [0; A] < [0; AA] < \dots < [0; A^k] < \dots < [0; A^{k+1}] < \dots$$
there exists $k \geq 0$ such that 
$$[0; A^k] \leq y < [0; A^{k+1}];$$
hence, $y < [0; A^{k+1}] \leq [0; A + y] = x$ and $M_{a, a-1, M}^{-1}(x-1) \in (0, x)$, so 
$M_{a, a-1, M}= M_{x, x-1, M}$. 

The case $x \leq a$ is similar: the only non-negative element of the orbit this time is
$$M_{a, a, N}^{-1}(x) = y \quad \textup{with }0 \leq y < \alpha^-$$
which, since $\alpha^- < x$, still implies $M_{a, a, N} = M_{x, x, N}$.
\end{proof}

\textbf{Proof of theorem \ref{main}}
Let $x \in I_a$, $a$ maximal. By proposition \ref{matchrat}, $M_{a, a, N}$ and $M_{a, a-1, M}$ are related by 
the identity  $(N,M)_{\rm alg}$.
Since by lemma \ref{sameorbit}, $M_{x, x, N} = M_{a, a, N}$ and $M_{x, x-1, M} = M_{a, a-1, M}$, 
the algebraic matching condition  $(N,M)_{\rm alg}$ holds also for $x$. 
\vskip 0.5 cm

In order to complete the proof of theorem \ref{second}, we are left with proving that the entropy 
is monotone on every maximal $I_a$:

\begin{proposition} \label{monot}
Let $I_a$ be a maximal quadratic interval, and let $N$ and $M$ be as in theorem \ref{main}: then the function
$\alpha \mapsto h(T_\alpha)$ is:
\begin{itemize}
\item[(i)] strictly increasing if $N < M$  
\item[(ii)] constant if $N = M$
\item[(iii)] strictly decreasing if $N >M$
\end{itemize}
on the whole interval $I_a$.
\end{proposition}

\noindent The proof is just an adaptation of the one given in \cite{NakadaNatsui} (see appendix): 
let us just remark that we are able to establish explicit bounds for the domain of validity 
of their entropy formula, which was previously just claimed to work locally. Moreover, 
$N$ and $M$ are now given in terms of the c.f. expansion of $a$, so it becomes immediate 
to establish which of the cases (i)-(ii)-(iii) holds in a neighbourhood of any given rational number.

\subsection{Period doubling}\label{2x}

Another feature observed in \cite{TiProCarMar}, (sect. 4.2.) was the production of infinite chains of 
adjacent matching intervals via {\it period doubling}; more formally, 
 
\begin{proposition}
Let $a$ be the pseudocenter of a maximal interval $I_a$, and write $a = [0; A^-] = [0; A^+]$ with $|A^-| \equiv 1 \mod 2$.
Then $a' := [0; A^- A^-]$ is the pseudocenter of a maximal interval.
\end{proposition}

The proposition follows immediately from lemma \ref{adiacent}, which will be proved in next section.
By applying the proposition repeatedly, one gets the 

\begin{corollary}
Let $I_a$ be a maximal (hence matching) interval. Then there is a countable chain of matching intervals
$$ \dots < I_{a_{n+1}} < I_{a_n} < \dots < I_{a_1} = I_a$$ 
such that $I_{a_n}$ and $I_{a_{n+1}}$ are adjacent, and $\lim_{n \rightarrow \infty} a_n := a_\infty >0$. 
\end{corollary}

Note that the proposition also gives a recursive algorithm to generate the c.f. expansion of the limit point $a_\infty$:
an explicit computation for the chain generated by $I_{1/2}$ is contained in \cite{TiProCarMar}, sect. 4.2.

\section{String techniques}\label{strtec}

This section contains the proofs of a few technical lemmata about the string ordering mentioned 
in the rest of the paper.

\subsection{String formalism}

To prove our results we shall need to fix some notation to manipulate the strings of partial quotients.

If $A$, $B$ are two finite strings composed with the alphabet  $\NN_+$ we denote

\begin{narrower}

\def\MyLen{1.5truecm}

\begin{description}

\item[\hbox to \MyLen{\hfill$A'$\hfill}] the \emph{twin string} of the finite string $A$ i.e.
the string such that the finite c.f.'s $[0,A]$ and $[0,A']$ 
represent the same rational number;

\item[\hbox to \MyLen{\hfill$AB$\hfill}] the concatenation of $A$ and $B$; $Ab$ will denote the concatenation of the finite string $A$ with the one-letter string $(b)$;

\item[\hbox to \MyLen{\hfill$A^n$\hfill}]
the concatenation of $n$ copies of
$A$  ($A^0$ is the empty string);

\item[\hbox to \MyLen{\hfill$\overline{A}$\hfill}]
means the endless
concatenation of $A$;

\item[\hbox to \MyLen{\hfill$|A|$\hfill}]
the length of $A$

\item[\hbox to \MyLen{\hfill$(A)_i^j$\hfill}] the substring of $A$
  going from the $i$-th figure to the $j$-th figure of $A$; 
to indicate $j$-th figure of the string $A$
  we shall usually write $(A)_j$ instead of $(A)_j^j$; 

\item[\hbox to \MyLen{\hfill$A \subseteq B$\hfill}]  means that $A$ is a {\em prefix}
of $B$, i.e. there exists $B_1$ such that $B = AB_1$.
\end{description}

\end{narrower}

We will be interested in the alternating lexicographic order structure on the
 space of finite or infinite strings as defined in section \ref{coding}. 
Note that, the set of finite strings $\mathcal{S}$ is a semigroup for the
operation of concatenation.  Associating a finite string $S$ to the
fractional map $( x\mapsto [0;S+x])$ yields a natural action of the
semigroup $\mathcal{S}$ on $\mathbb{R}_+$. Let us also recall that the
map $( x\mapsto [0;S+x])$ is increasing if $|S|$ is even and
decreasing if $|S|$ is odd, in particular odd convergents of any $x$
are greater than $x$ while even convergents are smaller. Moreover, if
$x:=[0;S,a+x']$ and $y:=[0;S,b+y']$ with $a>b \in \mathbb{N}_+$, $x', y' \in
[0,1[$, then $x>y$ if $|S|$ is even and $x<y$ if $|S|$ is odd.

In the following we shall need some effective criterion to compare infinite
periodic strings $\overline{S}, \ \overline{T}$: as soon as $|S|\neq |T|$ this becomes a
nontrivial task. The next section will deal this issue.
 
\subsection{String Lemma}

\begin{lemma} 
Let $S, \ T$ be two nonempty strings. Then the pair of infinite strings
$\overline{ S}$, $\overline{ T}$ is ordered in the same way as the pair $ST$, $TS$; namely 
 $$ ST \gtreqqless TS \iff \overline{S} \gtreqqless \overline{ T}.$$
\end{lemma}

\begin{proof}
If $ST = TS$ we can prove that there exists another string $P$ and integers $k,h \in \NN$ such that $S=P^k$, $T=P^h$, hence $\overline{S} = \overline{T}$.
In fact, we proceed by induction on $n:= \max\{S,T\}$.
For $n=1$ the claim is obviously true. Assume now we have proved this
claim for all pairs of strings of length strictly less than $n$, and
let $S,T$ be a pair of strings of maximal length $n$. We may assume
that $0< |T| < |S| \leq n$, the cases $|T|=0$ and $|T|=|S|$ being trivial. The hypothesis $TS=ST$ implies that $T$ is a prefix of $S$,
namely $S=TS_1$ therefore $TS=ST$ translates into $TS_1=S_1T$. Since
$\max\{|T|,|S_1|\}<|S|\leq n$ we use the inductive hypothesis to conclude that
$T=P^k$, $S_1=P^h$, and therefore $S=P^{h+k}$.  

If $ST \neq TS$, then $d := \min \{ j \in \mathbb{N} : (ST)_j \neq (TS)_j \} \leq s + t$. By lemma \ref{STn=TSn} with $n = d-1$ one has 
$$(ST)_1^d = (\overline{S})_1^d \qquad (TS)_1^d = (\overline{T})_1^d$$
hence the pair $(\overline{S}, \overline{T})$ is ordered in the same way as $(ST, TS)$.
\end{proof}

\begin{lemma}\label{STn=TSn}
Let $S, \ T$ be two nonempty strings, $s:=|S|$, $t:=|T|$, $n\in \NN$, $0 \leq n < s+t$. If  
$(ST)_1^n=(TS)_1^n$ then 
$$
\left\{
\begin{array}{cccc}
(\overline{S})_1^{n+1} & = & (ST)_1^{n+1} & \ \ \ (*)  \\
(\overline{T})_1^{n+1} & = & (TS)_1^{n+1} & \ \ \ (**)  
\end{array}
\right.
$$  
\end{lemma}

{\bf Proof of lemma \ref{STn=TSn}.}
We can assume $|T| \leq |S|$. We can split the proof in three cases, depending on the relation
between $n$ and the lengths data $t$ and $s$.  

{\bf Case 1: $0 \leq n <
  t$.} In this case both (*) and (**) trivially hold.  

{\bf Case 2: $
  n < s$, $kt\leq n < (k+1)t$ for some $k\geq 1$.}  Hypothesis (i)
implies that $T^k$ is a prefix of $S$, i.e. $S=T^kS_1$.  On the other
hand
\begin{enumerate}
\item[-] $\overline{S}$ coincides with $ST$  on the first $s$ figures $\stackrel{n<s}{\Longrightarrow} $ (*) holds;
\item[-] $\overline{T}$ coincides with $TS$  on the first $(k+1)t$ figures $\stackrel{n<(k+1)t}{\Longrightarrow} $ (**) holds;
\end{enumerate}

{\bf Case 3: $ s \leq n < s+t$.}  Hypothesis (i) implies that $S$
is a prefix of $T^{k}$ (with $k=\lceil{\frac{s}{t}}\rceil$),
i.e. $S=T^{k-1}T_0$, $T = T_0 T_1$.  Thus
$$
(\overline{S})_1^{s+t}  =  T^{k-1} T_0 T_0 T_1 = ST \qquad (\overline{T})_1^{s+t} = T^{k} T_0 = TS
$$
So (*) and (**) are again both verified.

The following remark will be useful further on
\begin{remark}\label{easy}
Let $T,\  S$ be two nonempty strings and set $a:=[0;ST]$, $b:=[0;S]$, $I_a:=]\alpha^-,\alpha^+[$ and $I_b:=]\beta^-,\beta^+[$. Then
\begin{enumerate}
\item[(i)] If $|S|$ is even then $b<a$ and $\beta^-<\alpha^-$;
\item[(ii)] If $|S|$ is odd and $T\neq(1)$, then $b>a$ and $\beta^+>\alpha^+$
\end{enumerate}
\end{remark}

\begin{lemma}
If $I_a \cap I_b \neq \emptyset$, $I_a \setminus I_b \neq \emptyset$ and $I_b \setminus I_a \neq \emptyset$, then 
either $I_a$ or $I_b$ is not maximal. 
\end{lemma}

\begin{proof}
By Lemma \ref{lemmacontained}, 
without loss of generality, we may assume that $a$ is
a convergent of $b$; hence we can write $a = [0; A]$, $b = [0; A^\ell A_0]$, where $A_0 \neq
\emptyset$ is a proper prefix of $A$.  Let $a_0:=[0;A_0]$; we claim
that the interval $I_{a_0}$ contains either $I_a$ or $I_b$.
There are several cases
to be examined; in all cases the proof
that the two intervals are nested, one inside the other, amounts to checking two inequalities: 
one of the two inequalities will be a trivial consequence of the previous remark 
while the other  is harder, but it will follow from the String Lemma \ref{STvsTS}.
We treat just one case in detail, and provide a table explains how to get the ``hard'' inequality for all the other cases.
Let $|A| \equiv 0$, $|A_0| \equiv 0$, $\alpha^+ = [0; \overline{A}]$, $\beta^+ = [0; \overline{A^\ell A_0}]$.
$$\alpha^+ < \beta^+ \Leftrightarrow \overline{A} < \overline{A^\ell
  A_0} \Leftrightarrow AA_0 < A_0 A \Leftrightarrow \overline{A} <
\\ \overline{A_0}$$ so $\alpha^+<\alpha_0^+$ and, by the remark
\ref{easy}, $\alpha_0^-<\alpha^-$ so that $I_{a} \subseteq I_{a_0}$
where $a_0 := [0; A_0]$.

\bigskip

\noindent
\!\!\!\!\!\!\!\!\!\!\!\!\!\!\!\!\!\!\!\!\!\!\!\!\!\!\!\!\begin{tabular}{|ll|c|ll|c|}
\hline
&&&&&\\
Cases & &  hypotheses used &  hard inequalitiy  && aim\\
\hline
&&&&&\\
$|A|$ even  & $|A_0|$ even & $\overline{A^\ell A_0}>\overline{A}$ & $\overline{A} <\overline{A_0}$ & & \\ 
$a<b$ & $a_0<a<b<\alpha^+$ & $\beta^+>\alpha^+ $ &$ \alpha^+<\alpha_0^+$ & & $I_a \subset I_{a_0}$\\ 
$\alpha^+:=[0;\overline{A}]$ & & & & &\\
\hline
&&&&&\\
$|A|$ even  & $|A_0|$ odd & $\overline{A^\ell A_0}<\overline{A}$ & $\overline{A_0} <\overline{A}$ & & \\ 
$a<b$ & $\alpha^+<b<a_0$ & $\beta^-<\alpha^+ $ &$ \alpha_0^-<\beta^-$ & & $I_b \subset I_{a_0}$\\ 
$\alpha^+:=[0;\overline{A}]$ & & & & &\\
\hline

&&&&&\\
$|A|$ odd  & $|A^\ell A_0|$ even & $\overline{A^\ell A_0}>\overline{A}$ & 
$\left\{
\begin{array}{l}
\overline{A_0} >\overline{A^\ell A_0}\\
\alpha_0^+>\beta^+
\end{array}
\right.$
 &if 
$\left\{
\begin{array}{cc}
|A_0| &\ {\rm even} \\
\ell  &\ {\rm even}
\end{array}
\right.$
 & $I_b \subset I_{a_0}$\\[15pt] 

$b<a$ & $b<\alpha^-$ & $\beta^+>\alpha^- $ 
&
$\left\{
\begin{array}{l}
\overline{A_0} <\overline{A}\\
\alpha_0^-<\alpha^-
\end{array}
\right.$
 &if 
$\left\{
\begin{array}{cc}
|A_0| &\ {\rm odd} \\
\ell  &\ {\rm odd}
\end{array}
\right.$
& 
$I_a \subset I_{a_0}$\\ 
$\alpha^-:=[0;\overline{A}]$ & & & & &\\
\hline

&&&&&\\
$|A|$ odd  & $|A^\ell A_0|$ odd & $\overline{A^\ell A_0'}>\overline{A}$ & 
$\left\{
\begin{array}{l}
\overline{A_0'} >\overline{A^\ell A_0'}\\
\alpha_0^+>\beta^+
\end{array}
\right.$
 &if 
$\left\{
\begin{array}{cc}
|A_0'| &\ {\rm even} \\
\ell  &\ {\rm odd}
\end{array}
\right.$
 & $I_b \subset I_{a_0}$\\[15pt] 

$b<a$ & $\alpha^-<b$ & $\beta^+>\alpha^- $ 
&
$\left\{
\begin{array}{l}
\overline{A_0'} <\overline{A}\\
\alpha_0^-<\alpha^-
\end{array}
\right.$
 &if 
$\left\{
\begin{array}{cc}
|A_0'| &\ {\rm odd} \\
\ell  &\ {\rm odd}
\end{array}
\right.$
& 
$I_a \subset I_{a_0}$\\ 
$\alpha^-:=[0;\overline{A}]$ & & & & &\\
\hline
                       
\hline                       
\end{tabular} 

\bigskip

\end{proof}

\begin{lemma}\label{adiacent}
Let $a_1=[0,P]$,  $a_\ell=[0,P^\ell]$.  The following are equivalent
\begin{enumerate}
\item[(i)] $I_{a_\ell}$ is maximal;
\item[(ii)]  $I_{a_1}$ is maximal and
$$\begin{array}{ll}
\ell=1 & \mbox{ if } |P| \mbox{ is even},\\
\ell \leq 2 & \mbox{ if } |P| \mbox{ is odd}.\\
\end{array}$$
\end{enumerate}
\end{lemma}

\begin{proof}
$[(i)\Rightarrow (ii)]$.
If $|P|$ even, $\ell>1$ and $a_{\ell-1}=[0;P^{\ell-1}]$, then
$I_{a_{\ell-1}}\supsetneq I_{a_{\ell}}$ so that $I_{a_\ell}$ can't be maximal.
If $|P|$ is odd and $\ell>2$, setting $a_{\ell-2}=[0;P^{\ell-2}]$ then
$I_{a_{\ell-2}}\supsetneq I_{a_{\ell}}$ so, again,  $I_{a_\ell}$ can't be maximal.
To conclude the proof we just need to prove that $I_{a_1}$ is maximal.
Let $I_{a_*}$ be the maximal interval containing $I_{a_1}$, so that $a:=[0;P_*]$ is a convergent of $a_1$.
The function $\phi(x):=[0;P+x]$ is injective, $\phi : I_{a_*}\stackrel{\sim}{\rightarrow}\phi(I_{a_*})= I_{\phi(a_*)}$, with
$\phi(a_*):=[0;PP_*]$; moreover $\phi(I_{a_1})= I_{\phi(a_1)}=I_{a_2}$. So
$$ I_{a_1}\subset I_{a_*}, \ \ \ \ \phi(I_{a_1})\subset \phi(I_{a_*})=I_{\phi(a_*)}, \ \ \ \ I_{a_2}=\phi(I_{a_1})\subset \phi(I_{a_*})=  I_{\phi(a_*)}.$$
Since $I_(a_2)$ is  maximal, $I_{a_2}=I_{\phi(a_*)} $ and thence $I_{a_1}=I_{a_*}$ is maximal.

$[(ii)\Rightarrow (i)]$. Let $|P|$ be odd and $\ell=2$ (otherwise there's nothing to prove!); we have to show that $I_{a_2}$ is maximal (if $I_{a_1}$ is).
Let $I_{a_j}:=]\alpha_j^-, \alpha_j^+[$ ($j=1,2$) and observe that, since $a_1$ is an odd convergent of $\alpha_1^-$ and $a_1$ is an even convergent of $\alpha_2^+$,
 $$a_2:=[0;PP]<[0;\overline{P}]=\alpha^+=\alpha^-<a_1.$$
If $I_a$ is the maximal interval containing $I_{a_2}$, $a:=[0;A]$, $|A|$ even,
  we have that $I_a\cap I_{a_1}=\emptyset$ and so $\overline{A}=\overline{P}$.
Therefore $A$ and $P$ have a common period $Q$: $A=Q^m$, $P=Q^\ell$; on the other hand, by virtue of the implication $[(i)\Rightarrow (ii)]$
 (the one we have already proved) we get $\ell=1$ ($\ell=1$ is impossible, since $|P|$ is odd) and therefore $m=2$, so $I_{a_2}=I_a$ is maximal.
\end{proof}

Let $S, \ T$ be two nonempty strings and
\begin{equation}\label{notation}
 \begin{array}{ccc}
a:=[0;ST], & b:=[0;S], & c:=[0;T];\\[12pt]
I_a:=]\alpha^-,  \alpha^+[, &  I_b:=]\beta^-, \beta^+[, & I_c:=]\gamma^-,\gamma^+[.
\end{array}
\end{equation}

\begin{proposition} \label{SLSAA}
Let $a = [0; A] \in \mathbb{Q} \cap (0, 1]$. 
\begin{enumerate} 
\item 
The following are equivalent:
\begin{itemize}
 \item[(i)] $I_a$ is maximal.
\item[(ii)] If $A = ST$ with $S$, $T$ finite nonempty strings, then either
$ST < TS$ or $ST = TS$ with $T = S$, $|S|$ odd .
\end{itemize}
\item  
Moreover, if $a=[0;ST]$ and $I_a$ is maximal, $[0; T] > [0; ST]$ (i.e. $a<c$).
\end{enumerate}
\end{proposition}

\begin{proof}
(1)--$[(i)\Rightarrow (ii)]$. Let us use the notation introduced above in
  \eqref{notation}; by maximality of $I_a$ we immediately get that $b
  \notin I_a$ ($a,b\in \mathbb{Q}$ and $den(b)< den(a)$ - see also definition \ref{pseudoc}); since $b\in I_b\setminus I_a \neq \emptyset$,
  maximality of $I_a$ and lemma \ref{strongmil} also imply that
  $I_a \cap I_b = \emptyset$.
\begin{enumerate}
\item[Case 0.]  If $b$ is an even convergent of $a$ (i.e. if $|S|$ is
  even and $b<a$) then $I_b$ lies to the left of $I_a$ and hence 
$\beta^+\leq  \alpha^-$; since $[0;\overline{S}]= \beta^+$ and
    $[0;\overline{ST}]\in \{\alpha^\pm\}$, by String Lemma \ref{STvsTS}
      we get $SST\leq STS$ and, since $|S|$ is even $ST\leq TS$. Lemma \ref{adiacent} tells us that,
since $|S|$ is even and $I_a$ is maximal, equality cannot hold.
 
\item[Case 1.] If $b$ is an odd convergent (i.e. if $|S|$ is odd and $b>a$) by the previous argument
  $\alpha^+\leq \beta^-$.  If
  $[0;\overline{ST}]=\alpha^+=\beta^-=[0;\overline{S}]$ then, by lemma
  \ref{adiacent}, $T=S$. If not, then
  $[0;\overline{ST}]<[0;\overline{S}]$; by String Lemma \ref{STvsTS},
  $STS<SST$ and, since $|S|$ is odd, this implies that $TS>ST$ (which
  is the same conclusion as the previous case).
\end{enumerate}
The first implication is thus proved.

(1)--$[(ii)\Rightarrow (i)]$.  Assume $I_a$ {\bf is not} maximal; then
there exist two non-empty strings such that $a:=[0;ST], \ b:=[0;S],
\ \ I_b$ is maximal and $I_b \supset I_a$ (which, in particular,
implies that if $|S|$ is odd then $S\neq T$). Then $\alpha^+\leq
\beta^+$ and $\alpha^-\geq\beta^-$. Let us give a quick glance at the
cases that can occur:

\bigskip

\begin{tabular}{|c|c|c|c|c|c|}
\hline
&&&&&\\
$|S|$ & $|T|$ & $[0;\overline{ST}]$ & $[0;\overline{S}]$ & consequence of String Lemma & Conclusion\\
&&&&&\\
\hline
&&&&&\\
even  &  even & $\alpha^+$ & $\beta^+$ & $STS \leq SST$ & $ST\geq TS$ \\ 
even  &  odd & $\alpha^-$ & $\beta^+$ &$STS < SST$ & $ST> TS$ \\ 
odd  &  even & $\alpha^-$ & $\beta^-$ & $STS \leq SST$ & $ ST\geq TS$ \\ 
odd  &  odd & $\alpha^+$ & $\beta^-$ & $STS < SST$ & $ ST > TS $ \\
\hline
\end{tabular} 

\medskip
It is thus easy to realize that condition (ii) never holds.
\bigskip

(2) Let us now prove the second statement of the previous proposition. Since our claim concerns rational values,
we may assume that $|ST|$ is even (so that
$\alpha^+=[0;\overline{ST}]$). Let us rule out the ``period doubling case'' (i.e. $|S|$ odd and $S=T$): in this case $a<c$ because $c$ is an odd convergent of $a$.
In all other cases the strict inequality  $ST<TS$ holds and hence $STT<TST$.

Moreover we know that
\begin{itemize}
\item
$\gamma:=[0;\overline{T}]$ is an endpoint of $I_c$; 
\item
$\gamma > \alpha^+$ (because $STT\leq TST$);
\item
$I_a\cap I_c =\emptyset$ because  $I_c$ must contain points which are not in $I_a$, and $I_a$ is maximal (recall lemma \ref{strongmil}).
\end{itemize}
Therefore $c>a$  (and in fact $\alpha^+\leq \gamma^-$ since $I_a\cap I_c=\emptyset$). 
\end{proof}

Let us point out that proposition \ref{SLSAA} provides an effective algorithm
to decide whether or not a string defines the pseudocenter of a
maximal interval: it is sufficient to check that all its cyclical
permutation produce strings which are strictly bigger (except if the
exceptional case of period doubling occurs).

\section*{Appendix}
{\bf (A) Comparison between matching conditions}
\medskip

Let us recall the matching conditions given in \cite{NakadaNatsui}:

\begin{itemize}
\item[\textup{(c-1)}] $\{T^n_{\alpha}(\alpha) : 0 \leq n < k_1 \} \cap \{T^m_{\alpha}(\alpha -1) : 0 \leq m < k_2 \} = \emptyset$
\item[\textup{(c-2)}] $M_{\alpha, \alpha, k_1} = \matr{1}{1}{0}{1} M_{\alpha, \alpha -1, k_2}$  $( \Rightarrow T^{k_1}_{\alpha}(\alpha) = T^{k_2}_{\alpha}(\alpha-1))$
\item[\textup{(c-3)}] $T^{k_1}_{\alpha}(\alpha)(= T^{k_2}_{\alpha}(\alpha -1)) \notin \{\alpha, \alpha -1\}$
\end{itemize}

\noindent The matching set is therefore
$$\mathcal{M}_{NN} := \{ \alpha \in (0, 1) : \textup{(c-1), (c-2), (c-3) hold for some }(k_1, k_2) \}$$

\begin{proposition}
If $\alpha$ satisfies the algebraic matching condition  $(N,M)_{\rm alg}$,
 then $T^{N+1}_{\alpha}(\alpha) = T^{M+1}_{\alpha}(\alpha-1)$.
\end{proposition}

\begin{proof}
By writing the identity  $(N,M)_{\rm alg}$ in terms of M\"obius transformations and evaluating it at $\alpha$, 
$$T_\alpha^N(\alpha) + T_\alpha^M(\alpha-1) =  -T^N_\alpha(\alpha) T^M_\alpha(\alpha-1)$$
which implies $T^N_\alpha(\alpha) = 0 \Leftrightarrow T^M_\alpha(\alpha-1) = 0$. If both are zero, 
the claim follows trivially since $T_\alpha(0) = 0$; if they are nonzero, one can write
\begin{equation}
 \label{matchinv}
 \frac{1}{T^N_\alpha(\alpha)}+ \frac{1}{T^M_\alpha(\alpha-1)} =- 1
\end{equation}
Now suppose $\epsilon_\alpha(T^N_\alpha(\alpha)) = \epsilon$, and $c_\alpha(T^N_\alpha(\alpha)) = c$
so that 
$$\frac{\epsilon}{T^N_\alpha(\alpha)} - c \in [\alpha-1, \alpha)$$
The fact that $|T^N_\alpha(\alpha)| < 1$ and \eqref{matchinv} imply $\epsilon_\alpha(T^M_\alpha(\alpha-1)) = -\epsilon$, 
hence 
$$-\frac{\epsilon}{T^M_\alpha(\alpha-1)} - c -\epsilon \in [\alpha-1, \alpha)$$
so $c_\alpha(T^M_\alpha(\alpha-1)) = c + \epsilon$ and $T^{N+1}_\alpha(\alpha) = T^{M+1}_\alpha(\alpha-1)$.
\end{proof}

\begin{proposition}
Let $I_a$ be a maximal quadratic interval, and let the two c.f. expansions of $a$ be 
$a = [0; A^+] = [0; A^-]$.
Let $N$ and $M$ as is theorem \ref{main} and
$$\tilde{I}_a := \{ \alpha \in I_a \textrm{ s.t. (c-1), (c-2), (c-3) hold with }k_1 = N+1, k_2 = M+1 \}$$
Then 
$$I_a \setminus \tilde{I}_a \subseteq \{a \} \cup \{ \alpha = [0; \overline{A^+, k}], k \in \mathbb{N} \} \cup \{ \alpha = [0; \overline{A^-, k}], k \in \mathbb{N} \}$$
\end{proposition}

\begin{proof}
By the proof of the previuos proposition, (c-2) holds for $\alpha \in I_a \setminus \{a\}$. 
By using the explicit description of the orbits as in corollary \ref{orbitpseudo} and lemma \ref{sameorbit}, one can check (c-1) holds for every $\alpha \in I_a \setminus\{a\}$.
Exceptions to (c-3) precisely correspond to $\alpha = [0; \overline{A^-, k}]$ or $\alpha = [0; \overline{A^+, k}]$. 
\end{proof}

\begin{corollary}
 $\mathcal{M} \setminus \mathcal{M}_{NN}$ is a countable set.
\end{corollary}

\bigskip
\noindent {\bf (B) The entropy is monotone on maximal intervals}
\medskip 

Let us now prove proposition \ref{monot}:
\begin{lemma}
Let $a \in \mathbb{Q} \cap (0, 1]$ such that $I_a = (\alpha^-, \alpha^+)$ is maximal, let 
$a = [0; A]$ be its c.f. expansion with $|A| \equiv 0 \mod 2$, and choose $\alpha$, $\alpha'$
such that $\alpha^- < \alpha < \alpha' < \alpha^+$ and $\alpha' \leq [0; A + \alpha]$.
Then  
$$\frac{h(T_{\alpha'})}{h(T_\alpha)} = 1 + (M-N)\mu_{\alpha'}([\alpha, \alpha']) $$
$$\frac{h(T_{\alpha})}{h(T_{\alpha'})} = 1 + (N-M) \mu_{\alpha}([\alpha-1, \alpha'-1])$$ 
where $\mu_\alpha$ and $\mu_{\alpha'}$  are  the invariant densities of $T_\alpha$ and $T_{\alpha'}$, respectively.
\end{lemma}

\begin{proof} 
Choose $x \in (\alpha, \alpha')$. The proof proceeds exactly as in \cite{NakadaNatsui}, thm. 2, once we show that
$$M_{\alpha', x, k}^{-1}(x) \notin (\alpha, \alpha') \qquad 1 \leq k \leq N$$
$$M_{\alpha, x-1, h}^{-1}(x-1) \notin (\alpha-1, \alpha'-1) \qquad 1 \leq h \leq M$$

This follows directly from lemma \ref{anatomy}, except for two cases:
one in which $h = M$ and $x \geq a$, and the other in which and $k =
N$ and $x \leq a$. In the first case one can write $x = [0; A + y]$
with $a = [0; A]$, $|A| \equiv 0 \mod 2$, $0 \leq y < \alpha^+$. Then
$M_{\alpha, x-1, M}^{-1}(x-1) = y < \alpha$ because
$$[0; A + y] = x < \alpha' \leq [0; A + \alpha] \Rightarrow y < \alpha$$  
The second case is handled similarly.
\end{proof}

\textbf{Proof of proposition \ref{monot}}
Given $\alpha, \alpha' \in I_a$, $\alpha < \alpha'$, let 
$\alpha_k := [0; A^k + \alpha]$ and $k_0 := \max \{k > 0 \textup{ s.t. } \alpha_k < \alpha' \}$.
One can apply the lemma to each consecutive pair of the chain $\alpha < \alpha_1 < \dots < \alpha_{k_0} < \alpha$.

\end{document}